\documentclass [11pt]{article}
\usepackage{amssymb,amsmath,comment}
\usepackage[dvips]{graphicx}

\newtheorem{theorem}{Theorem}[section]

\newtheorem{lemma}[theorem]{Lemma}

\newtheorem{remark}[theorem]{Remark}

\begin{document}
\title{Combinatorial properties of Circulant Hadamard matrices}
\author{Reinhardt Euler\\
Lab-STICC UMR CNRS 6285,
University of Brest,\\
20, Av. Le Gorgeu, C.S. 93837,
29238 Brest Cedex 3, France,\\
\\
Luis H. Gallardo, Olivier Rahavandrainy\\
Laboratoire de Math\'ematiques UMR CNRS 6205, University of Brest,\\
6, Avenue Le Gorgeu, C.S. 93837, 29238 Brest Cedex 3, France.\\
\\
e-mail: \{Reinhardt.Euler,Luis.Gallardo,Olivier.Rahavandrainy\}\\
@univ-brest.fr}
\date{}
\maketitle

\\
\\
\\
{\bf{Abstract}} We describe combinatorial properties of the defining row of a circulant Hadamard matrix by exploiting
its orthogonality to subsequent rows, and show how to exclude several particular forms of these matrices.
\newpage

\section{Introduction and basic definitions}
Throughout  this paper, the $i$-row of a matrix $H$ of size  $n \times n$ is denoted by $H_i$, $1 \leq i \leq n$.
A $(\pm 1)$-matrix $H$ of size $n \times n$ is called \emph{Hadamard matrix} if its rows are mutually orthogonal.
It is a \emph{circulant matrix} if for $H_1 = (h_1, \ldots, h_n)$, its $i$-row is given by
$H_i = (h_{1-i+1},\ldots, h_{n-i+1})$, the subscripts being taken modulo~$n$. In this case, the first row $H_1$
will be called the \emph{defining row} of $H$.

An example of a $4 \times 4$ circulant Hadamard matrix is the following:

\[
H = \left(\begin{array}{rrrr}
-1&1&1&1\\
1&-1&1&1\\
1&1&-1&1\\
1&1&1&-1
\end{array} \right)
\]

A conjecture of Ryser (cf. \cite[pp.97]{Davis}) states that there is no circulant Hadamard matrix unless $n = 1$ or $n= 4$.
This conjecture is still open although many partial results have been obtained. For recent results and further details,
we refer to \cite{Craigen}, \cite{EGR} and \cite{Leung}.\\

The objective of this paper is to exhibit some combinatorial properties of the defining row $H_1$ of such a matrix $H$
that we obtained by exploiting the orthogonality of $H_1$ to some of its subsequent rows. It turns out that these
properties allow to exclude several particular forms for circulant Hadamard matrices. We think that a further development of this ``combinatorial approach'' may contribute to exclude other cases of such matrices and help to
determine the maximum number $k$ for which the first $k$ rows of a circulant $(\pm 1)$-matrix are mutually orthogonal.\\

The order $n$ (for $n > 1$) of a circulant Hadamard matrix $H$ is well known to be of the form $n = 4m=4h^2$, where
$h$ is an odd integer. Also, by orthogonality, the entries of any two different rows have the same sign in exactly
$2m$ columns, and a different sign in the other $2m$ ones.\\

Finally, the property of being circulant allows any row of $H$ to be represented as a circular sequence of ``$+$''
and ``$-$''.

Some further definitions are required. Let us consider the first row $H_1$ of a circulant Hadamard matrix $H$.
A maximal sequence of consecutive $+1$'s (respectively $-1$'s) will be called a $\emph{positive block}$
(respectively $\emph{negative block}$), denoted by $B$. If the size of such a block equals $k$, for some positive
integer $k$, we will speak of a {\em positive} (respectively {\em negative}) {\em $k$-block}. The family of all blocks is denoted by
${\cal{B}}$, that of all $k$-blocks by ${\cal{B}}_k$ and that of all $k$-blocks with $k \geq i$, for a given positive
integer $i$, by ${\cal{B}}_{\geq i}$. Thus, $H_1$ decomposes into a sequence of alternating, positive or negative blocks.
For convenience, we suppose the first and last block of $H_1$ to be of opposite sign. Finally, a maximal alternating
sequence of blocks in ${\cal{B}}_k$ (in ${\cal{B}}_{\geq k}$) will be called a {\em $k$-alternating sequence}
({\em $\geq k$-alternating sequence}, respectively), the {\em length} of which is given by the number of blocks forming it.\\

\section{Our results}

\begin{lemma} \label{blockline1}
The total number of blocks in the first row $H_1$ of a circulant Hadamard matrix $H$ equals $2m$.
\end{lemma}

{\bf{Proof}}: We consider $H_1$ and $H_2$ of $H$ and observe that the number of columns with entries of different sign
equals the total number of blocks. Hence, $|{\cal{B}}| = 2m$. \hfill{$\Box$} \\

Alternatively, if we count the number of columns whose entries have equal sign, we obtain
$\displaystyle{\sum_{B \in {\cal{B}}} (|B| -1)}$, which equals $4m- |{\cal{B}}|$ and which by orthogonality,
gives $2m$. Again, $|{\cal{B}}| = 2m$. \hfill{$\Box$} \\

\begin{lemma} \label{blockline2}
The total number of $1$-blocks in $H_1$ equals $m$.
\end{lemma}

{\bf{Proof}}: We consider $H_1$ and $H_3$ of $H$ and count the number of columns with entries of equal sign.
We observe that: \\

\noindent $\bullet\ $ in a $1$-alternating sequence, this number equals the length of the sequence minus $1$, plus $1$ arising from
the $\geq 2$-block preceding the sequence,\\
$\bullet\ $ a block $B \in {\cal{B}}_{\geq 3}$ contributes with $|B| - 2$, \\
$\bullet\ $ no block $B \in {\cal{B}}_{2}$ contributes to this number (except those preceding a $1$-alternating path,
whose contribution is already taken into account).\\

By orthogonality of $H_1$ and $H_3$, we obtain:

$$\begin{array}{lcl}
2m &= &\displaystyle{|{\cal{B}}_{1}| + \sum_{B \in {\cal{B}}_{\geq 3}} (|B| -2)}\\
&= &\displaystyle{|{\cal{B}}_{1}| + \sum_{B \in {\cal{B}}_{\geq 2}} (|B| -2)}\\
 &=& \displaystyle{|{\cal{B}}_{1}| +  \left(\sum_{B \in {\cal{B}}_{\geq 2}} |B|\right) - 2 |{\cal{B}}_{\geq 2}|} \\
 &=& \displaystyle{|{\cal{B}}_{1}| + \left(4m- |{\cal{B}}_{1}|\right) - 2(|{\cal{B}}| - |{\cal{B}}_{1}|)}\\
 &=& 4m-4m+2|{\cal{B}}_{1}|.
 \end{array}$$

\indent Hence, $|{\cal{B}}_{1}| = m$.\\

Alternatively, we may count the number of columns of $H_1$ and $H_3$, with entries of different signs. We denote by
$\alpha_1$ the number of $1$-alternating sequences, by $\alpha_{\geq 2}$ the number of ${\geq 2}$-alternating sequences,
by $\Lambda$ the set of all ${\geq 2}$-alternating sequences, and by $|\ell|$ the length of an element $\ell$ of $\Lambda$.
We observe that: \\

\noindent $\bullet\ $ any 1-alternating sequence contributes with 2,\\
$\bullet\ $ any ${\geq 2}$-alternating sequence $\ell$ contributes with $2(|\ell| - 1)$. \\

Thus, we obtain:

\begin{equation}\label{equ1}
2m = 2 \alpha_1 + 2 (\sum_{\text{$\ell \in \Lambda$}} (|\ell| - 1)).
\end{equation}

Since 1-alternating sequences and ${\geq 2}$-alternating sequences are themselves alternating within $H_{1}$,
$\alpha_1 = \alpha_{\geq 2}$, and we get from (\ref{equ1}):

$$m = \alpha_1 + |{\cal{B}}_{\geq 2}| - \alpha_{\geq 2} = |{\cal{B}}_{\geq 2}| = |{\cal{B}}| - |{\cal{B}}_{1}| =
2m - |{\cal{B}}_{1}| \text{ by Lemma \ref{blockline1}}.$$

\indent Hence, $|{\cal{B}}_{1}| = m$. \hfill{$\Box$}

\begin{lemma} \label{blockline3}
The number of blocks of size $2$ and the number of $1$-alternating sequences add up to $m$:
$$|{\cal{B}}_2|+ \alpha_1 = m.$$
\end{lemma}

{\bf{Proof}}:  We consider $H_1$ and $H_4$ of $H$ and count the number of columns with entries of equal sign.
We observe that: \\

\noindent $\bullet\ $ any $1$- alternating sequence is preceded and followed by a block of size~$\geq 2$,\\
$\bullet\ $ a block $B \in {\cal{B}}_{\geq 3}$ contributes with $|B| - 3$, \\
$\bullet\ $ any $2$-block $B \in {\cal{B}}_{2}$ is preceded and followed by an element of opposite sign.\\

By orthogonality of $H_1$ and $H_4$ and by Lemmata \ref{blockline1} and \ref{blockline2}, we obtain:
$$\begin{array}{lcl}
2m &= &\displaystyle{\sum_{B \in {\cal{B}}_{\geq 3}} (|B| -3)} + |{\cal{B}}_{2}| + 2 \alpha_{1}\\
 &=& 4m-\displaystyle{|{\cal{B}}_{1}| - 2 |{\cal{B}}_{2}|} -  3 |{\cal{B}}_{\geq 3}|  + |{\cal{B}}_{2}| + 2 \alpha_{1} \\
 &=& 3m - |{\cal{B}}_{2}| - 3(|{\cal{B}}| -m- |{\cal{B}}_{2}|) + 2 \alpha_{1} \\
 &=& 2|{\cal{B}}_{2}| + 2 \alpha_{1}.
 \end{array}$$
\indent Hence, $|{\cal{B}}_{2}| + \alpha_{1}= m$. \hfill{$\Box$}

\begin{remark}
It would be interesting to know whether counting the number of columns with entries of different sign provides
the same condition on $|{\cal{B}}_2|$, similar to that obtained for $|{\cal{B}}_1|$. For the sake of completeness and
to make the reader more familiar with our approach, we have included such a condition at this place.
\end{remark}

For given $n$, let $\alpha_{1}, \alpha_{2}$ and $\alpha_{\geq 3}$ be the number of $1$-, $2$- and
$\geq3$-alternating sequences, and let $\alpha_{2, \geq 3}$ denote the number
of pairs of blocks $B, B'$ with $B \in {\cal{B}}_{2}$, $B' \in {\cal{B}}_{\geq 3}$ or vice versa.
We obtain the following:

\pagebreak

\begin{lemma} \label{blockline4}
$$\alpha_1 + \alpha_{2, \geq 3}= \alpha_2 + \alpha_{\geq 3}.$$
\end{lemma}

{\bf{Proof}}:  Any $1$-, $2$- or $\geq 3$-alternating sequence is preceded (followed) by a $2$- or $3$-, $1$-
or $\geq 3$-, $1$- or $2$-alternating sequence, respectively.\\
Let us examine the first case, i.e. that of a $1$-alternating sequence. Counting the number of columns with
entries of different signs gives us a contribution of $(\ell - 1)$ for such a
sequence, if $\ell$ denotes the length of it. \\

For a $2$-alternating sequence preceded and followed by a $1$-alternating one, the total contribution
will be $(\ell - 1)$ again. If, however, such a sequence is preceded or followed by a $\geq 3$-alternating
sequence (see Figure 3), one or two consecutive pairs
of blocks $B, B'$ with $B \in {\cal{B}}_{2}$, $B' \in {\cal{B}}_{\geq 3}$ or vice versa, will appear and contribute
by an additional value of $1$. \\

In that case, there is another contribution of $1$, that will be included in the
contribution of the $\geq 3$-alternating sequence. \\

Finally, for $\geq 3$-alternating sequences: The overall contribution
of such a $\geq 3$-alternating sequence is $3\ell - 1$. Also, observe that the pair of
blocks $B, B'$ with $B \in {\cal{B}}_{\geq 3}$, the last block of this sequence, and $B' \in {\cal{B}}_{2}$,
the first of the following $2$-alternating one, again contribute by $1$ to the value $\alpha_{2,\geq 3}$. \\

Altogether, counting the number of columns with entries of different signs in $H_1$ and $H_4$ gives the condition:
$$\sum_{i \in I} (\ell_i - 1) + \sum_{j \in J} (\ell_j - 1) +  \sum_{k \in K} (3\ell_k - 1) + \alpha_{2,\geq 3} = 2m,$$

where $I, J, K,$ denote the sets of 1-, 2- and $\geq 3-$alternating sequences, respectively, and $\ell$ their length. \\

We obtain: \\

$$|{\cal{B}}_1| - \alpha_1 + |{\cal{B}}_2| - \alpha_2 + 3 |{\cal{B}}_{\geq 3}| - \alpha_{\geq 3} +  \alpha_{2,\geq 3} =
2m.$$

Since, by Lemma \ref{blockline1},
$$|{\cal{B}}_{1}|+|{\cal{B}}_{2}|+3|{\cal{B}}_{\geq 3}|=4m-2|{\cal{B}}_{1}|-2|{\cal{B}}_{2}|,$$

$$4m-\alpha_1-\alpha_2 - \alpha_{\geq 3} +  \alpha_{2,\geq 3} = 2|{\cal{B}}_1| + 2 |{\cal{B}}_2|.$$

By Lemma \ref{blockline2}, we obtain:

$$2m -\alpha_{1}-\alpha_{2} - \alpha_{\geq 3} +  \alpha_{2,\geq 3} = 2 |{\cal{B}}_{2}|.$$

Since, by Lemma \ref{blockline3},
$$2m - 2 |{\cal{B}}_{2}| = 2 \alpha_{1},$$
we obtain as stated
$$\alpha_1 + \alpha_{2, \geq 3}= \alpha_2 + \alpha_{\geq 3}.$$
\hfill{$\Box$}

\vspace{.8cm}

We are now ready to show the nonexistence of a circulant Hadamard matrix for the following four particular situations:
$$\text{i) $\alpha_ 1 = 1$, ii) $\alpha_1 = m$, iii) $\alpha_1 = m-1$, iv) $\alpha_1 = 2$}$$

\begin{theorem} \label{mainresult}
If in $H_1$, there is only one $1$-alternating sequence (of length $m$), or if there are $m$ such sequences (all of length $1$), or if the number
of $1$-alternating sequences is $m-1$, then there is no circulant Hadamard matrix $H$ with $H_1$ as first row.
\end{theorem}

{\bf{Proof}}: We consider the three cases.\\

{\underline{Case $\alpha_1 = 1$}}:\\

\noindent We know already, from Lemmata \ref{blockline1}, \ref{blockline2} and \ref{blockline3}, that:\\
- the total number of blocks is $2m$,\\
- that of $1$-blocks is $m$,\\
- $|{\cal{B}}_2| + \alpha_1 = m$.\\
Therefore, since $n=4m$, there is exactly one more block of size $m+2$, and if we suppose that this particular block precedes the
$1$-alternating sequence (clockwise), we obtain a first representation of $H_1$ for $\alpha_1 = 1$\\

\noindent Evaluation of the scalar product of $H_1$ and $H_5$ gives:
$$(m-4) + 0 + 2(m-3) + 0 + (m-2)+0 = 4m-12,$$
which allows for orthogonality only if $m=3$, contradicting the condition on $m$ to be an odd square integer.\\

By symmetry, we obtain the same result if the particular $(m+2)$-block follows the $1$-alternating sequence.\\

If, now, the $(m+2)$-block is situated ``in between'' the sequence of $2$-blocks, \\

\noindent we may suppose that $$m-1 = k_1 + k_2, \text{with $k_1, k_2$ both even or both odd}.$$
Evaluating the scalar product of $H_1$ and $H_5$ gives:
$$(m-4) + 0 + 2(k_1-2) + 0 + (m-2) + 2(k_2-2) + 0 = 4m-16.$$
Orthogonality implies $m = 4$, contradicting the condition on $m$ to be odd.\\
We obtain the same result for the case $k_1 = 1$ or $k_2 = 1$, so that altogether, $H_1$ cannot be the defining row of a circulant Hadamard
matrix.\\

{\underline{Case $\alpha_1 = m$}}:\\

By Lemma \ref{blockline3}, there is no $2$-block in $H_1$ and $m$ blocks remain to cover $3m$ elements, i.e. besides the $m$ $1$-blocks, there are
exactly $m$ $3$-blocks, and both types alternate in $H_1$. But then the scalar product of $H_1$ and $H_5$ gives $4m > 0$ and $H_1$ cannot define a circulant Hadamard matrix either.\\

{\underline{Case $\alpha_1 = m-1$}}:\\

In this case there is exactly one consecutive pair of $1$-blocks, and by Lemma \ref{blockline3}, there is exactly one $2$-block, so that $3m-2$ elements have to be covered by $m-1$ blocks of size $\geq 3$. This is only possible with $(m-2)$ $3$-blocks and one $4$-block. Moreover, the unique $2$-block, say $B$, cannot be preceded {\em and} followed by a $\geq 3$-block. Therefore, we come up with two possibilities:\\

a) $B$ is preceded and followed by a $1$-alternating sequence.\\

b) $B$ is preceded or followed (but not both) by a $\geq 3$-block.\\

Just observe, that in situation a) exactly two $\geq 3$-blocks follow each other. We count the number of columns
with different signs in $H_{1}$ and $H_{5}$ and observe that a contribution (of 2) to this number is given only by

\begin{itemize}
\item the two consecutive $1$-blocks;
\item the unique $2$-block;
\item the unique $4$-block;
\item exactly one of the $3$-blocks.
\end{itemize}

\noindent By orthogonality, $$2m = 8,$$
a contradiction to the condition on $m$ to be odd.\\

In situation (b), no $3$-block contributes any more, and orthogonality of $H_{1}$ and $H_{5}$ gives
$$2m = 6,$$ a contradiction to the condition that $m$ has to be an odd square integer.    \hfill{$\Box$}  \\

In the following, we will show how the orthogonality of $H_1$ and $H_5$ can be exploited to obtain a further condition involving $|{\cal{B}}_2|$
and $|{\cal{B}}_3|$, which we will use afterwards to exclude the situation where $\alpha_1 = 2$.

For this, we need to describe the possible types of $1$-alternating and $2$-alternating sequences  (up to symmetry) \\

It is not difficult to verify that if we count the columns of $H_1$ and $H_5$ with equal sign, the contribution of these $6$ types will be $3, 2,
1$ and $\ell$, $\ell-1, \ \ell-2$, respectively. \\

For $2$-alternating sequences, we obtain $3$ similar types, whose contributions will be $2\ell$, $2\ell-1,
\ 2\ell-2$, respectively: \\

\noindent with the properties that

\begin{itemize}

\item type 1 is preceded and followed by a $\geq 3$-block;
\item type 2 is preceded by a $1$-alternating sequence and followed by a $\geq 3$-block (or vice versa);
\item type 3 is preceded and followed by a $1$-alternating sequence.
\end{itemize}

If now $\alpha_1^i$ denotes the number of $1$-alternating sequences of type $i$, for $i = 1, \ldots, 6$, and $\alpha_2^j$ the number
of $2$-alternating sequences of type $j$, for $j = 1, \ldots, 3$, we get \\

\begin{lemma} \label{blockline5}
If $H$ is a circulant Hadamard matrix, then one has
$$2 |{\cal{B}}_2| + 2 |{\cal{B}}_3| = (\alpha_1^2 + 2 \alpha_1^3) + (2\alpha_1^4 + 3 \alpha_1^5 + 4 \alpha_1^6) + (\alpha_2^2 + 2 \alpha_2^3).$$
\end{lemma}

{\bf{Proof}}:\\
We can describe the total number of columns of $H_1$ and $H_5$, whose entries have equal sign, by the following expression
\[
L:= 3\alpha_1^1 + 2 \alpha_1^2 + \alpha_1^3 + C + D,
\]
where
\[
C:= \sum_{i \in I} \ell_i + \sum_{j \in J} (\ell_j - 1) + \sum_{k \in K} (\ell_k - 2),
\]
and
\[
D:= \sum_{p \in P} 2\ell_p + \sum_{q \in Q} (2\ell_q-1) + \sum_{r \in R} (2\ell_r-2) + \sum_{B \in {\cal{B}}_{\geq 4}} (|B|-4) +
|{\cal{B}}_3|,
\]

\noindent expression in which $\displaystyle{\sum_{j \in J} (\ell_j - 1)}$, for instance, relates to all $1$-alternating sequences $S_j$, $j \in J$, whose contribution equals $\ell_j - 1$, i.e., which are of~type~$5$.\\

We know the following:

$$\begin{array}{l}
i) \ \displaystyle{\alpha_1^1 + \alpha_1^2 + \alpha_1^3 + \sum_{i \in I} \ell_i + \sum_{j \in J} \ell_j + \sum_{k \in K} \ell_k = |{\cal{B}}_1| =
m}, \text{ by Lemma \ref{blockline2}},\\
\\
ii) \ \displaystyle{\sum_{p \in P} \ell_p + \sum_{q \in Q} \ell_q + \sum_{r \in R} \ell_r =  |{\cal{B}}_2|},  \text{\hspace{.6cm}and} \\

iii) \ \displaystyle{\alpha_1^1 + \alpha_1^2 + \alpha_1^3 + \alpha_1^4 + \alpha_1^5 + \alpha_1^6 = \alpha_1}.
\end{array}$$

Therefore,

$$\begin{array}{lcl}
\displaystyle{3\alpha_1^1 + 2\alpha_1^2 + \alpha_1^3 + C}&=& 2\alpha_1^1 + \alpha_1^2 + |{\cal{B}}_1| - \alpha_1^5 - 2 \alpha_1^6\\
&=& \alpha_1 + m + \alpha_1^1 - \alpha_1^3 - \alpha_1^4 - 2 \alpha_1^5 - 3 \alpha_1^6.
\end{array}$$

We also have

$$\sum_{p \in P} 2\ell_p + \sum_{q \in Q} (2\ell_q-1) + \sum_{r \in R} (2\ell_r-2) = 2 |{\cal{B}}_2| - \alpha_2^2 - 2 \alpha_2^3.$$

Thus, if we put

$$\delta:= \alpha_1 + m + \alpha_1^1 - \alpha_1^3 - \alpha_1^4 -2 \alpha_1^5 - 3 \alpha_1^6 + 2 |{\cal{B}}_2| - \alpha_2^2 - 2 \alpha_2^3,$$

then

$$\begin{array}{lcl}
L &=& \displaystyle{\delta + \sum_{B \in {\cal{B}}_{\geq 4}} |B|- 4 |{\cal{B}}_{\geq 4}| + |{\cal{B}}_3|}\\
&=& \delta + 4m - |{\cal{B}}_1|- 2 |{\cal{B}}_2| - 3 |{\cal{B}}_3| - 4(|{\cal{B}}|-|{\cal{B}}_1|- |{\cal{B}}_2| - |{\cal{B}}_3|) +
|{\cal{B}}_3|\\
&=& \delta -m + 2|{\cal{B}}_2|+ 2 |{\cal{B}}_3|, \text{\ by Lemmata \ref{blockline1} and \ref{blockline2}}\\
&=& \delta - \alpha_1 + |{\cal{B}}_2|+ 2 |{\cal{B}}_3|, \text{\ by Lemma \ref{blockline3}}\\
&=& m + (\alpha_1^1 - \alpha_1^3 - \alpha_1^4 -2 \alpha_1^5 - 3 \alpha_1^6) + 3 |{\cal{B}}_2| + 2 |{\cal{B}}_3| - \alpha_2^2 - 2 \alpha_2^3.
\end{array}$$

This last expression, by orthogonality, has to be equal to $2m$.\\

So, by Lemma \ref{blockline3} and by iii) above, we get

$$(\alpha_1^1 - \alpha_1^3 - \alpha_1^4 -2 \alpha_1^5 - 3 \alpha_1^6)=(\alpha_1 - \alpha_1^2 - 2 \alpha_1^3 - 2 \alpha_1^4 - 3 \alpha_1^5 - 4 \alpha_1^6)$$
and

$$|{\cal{B}}_2| + \alpha_1 = m = \alpha_1 - \alpha_1^2 - 2 \alpha_1^3 - 2 \alpha_1^4 - 3 \alpha_1^5 - 4 \alpha_1^6 + 3 |{\cal{B}}_2| + 2
|{\cal{B}}_3| - \alpha_2^2 - 2 \alpha_2^3,$$

which gives our result:

$$2 |{\cal{B}}_2| + 2 |{\cal{B}}_3| = (\alpha_1^2 + 2 \alpha_1^3) + (2\alpha_1^4 + 3 \alpha_1^5 + 4 \alpha_1^6) + (\alpha_2^2 + 2 \alpha_2^3).$$
\hfill{$\Box$}\\
\\
Lemma \ref{blockline5} can now be used to exclude another situation:

\begin{theorem} \label{mainresult2}
If $H_1$ contains exactly two $1$-alternating sequences, i.e. if $\alpha_1 = 2$, then $H$ is not a circulant Hadamard matrix.
\end{theorem}

{\bf{Proof}}:\\
Lemma \ref{blockline3} implies that $|{\cal{B}}_2| = m-\alpha_1 = m-2$. There are two blocks left, by Lemma \ref{blockline1}, to cover the
remaining $m+4$ elements. Hence  $|{\cal{B}}_3| \in \{0,1\}$ for $m \geq 3$. Also, $\alpha_2^1 + \alpha_2^2 + \alpha_2^3 \leq 4$, with equality if
the two $1$-alternating sequences and the two remaining blocks from ${\cal{B}}_{\geq 3}$ are all alternating with a $2$-alternating sequence.

Since
$$(\alpha_1^2 + 2 \alpha_1^3) + (2\alpha_1^4 + 3 \alpha_1^5 + 4 \alpha_1^6) + (\alpha_2^2 + 2 \alpha_2^3) \leq 4 + 4 + 8 = 16,$$
\indent we get by Lemma \ref{blockline5}
$$m-2 + |{\cal{B}}_3| = |{\cal{B}}_2| + |{\cal{B}}_3| \leq 8, \text{ with $|{\cal{B}}_3| \in \{0,1\}$}.$$
Hence $m \leq 10$ or $m\leq 9$, both cases for which no circulant Hadamard matrix exists. \hfill{$\Box$}

\newpage

\section{Conclusion and further work}
Future work should include an analysis of further rows of $H$ along the same line. This is what we are currently doing.

\end{document}